\documentclass[reqno,11pt]{amsart}
\usepackage{amsmath,amsthm,amssymb}
\usepackage{latexsym,epsfig,lineno,hyperref}

\newtheorem{theorem}{Theorem}[section]
\newtheorem{lemma}[theorem]{Lemma}

\newtheorem{coro}{Corollary}[section]
\theoremstyle{definition}

\theoremstyle{remark}

\numberwithin{equation}{section}
\newcommand{\rb}{\right}
\newcommand{\lb}{\left}

\begin{document}

\title[Arithmetic properties of $2^\alpha-$Regular overpartition pairs]{Arithmetic properties of $2^\alpha-$Regular overpartition pairs}

\author{B. Hemanthkumar}
\address{Department of Mathematics, RV College of Engineering, RV Vidyanikethan Post, Mysuru Road, Bengaluru-560 059, Karnataka, India.}
\email{hemanthkumarb.30@gmail.com}
\author{ H. S. Sumanth Bharadwaj}
\address{Department of Mathematics and Statistics, M. S. Ramaiah University of Applied Sciences, RTC Campus, Peenya, Bengaluru - 560 058, Karnataka, India}
\email{sumanthbharadwaj@gmail.com}

\date{}

\begin{abstract}
Recently, several mathematicians have investigated various partition functions with the goal of discovering Ramanujan-type congruences. One such function is $\overline{B}_{2^\alpha}(n)$, which represents the number of $2^\alpha-$regular overpartition pairs of $n$. In this context, we establish Ramanujan-type congruences modulo powers of $2$ for this function. For instance, we prove that 
\begin{equation*}
\overline{B}_{2^{\alpha}}(2^{\alpha+\beta+1}(n+1)) \equiv 0\pmod{2^{3\beta+5}}
\end{equation*}
for all $n, \beta\geq 0,\, \alpha \in \mathbb{N}$.
\end{abstract}

\maketitle

\section{Introduction}\label{S1}
An \emph{overpartition} of a positive integer $n$ is a partition of $n$ 
in which the \emph{last} occurrence of a part may be overlined 
(see Corteel and Lovejoy \cite{CL}). 
For example, the partition
\[
   5 + \overline{5} + 3 + \overline{1}
\]
is an overpartition of 14, where the last occurrences of 5 and 1 are overlined. The number of overpartitions of $n$ is denoted by $\overline{p}(n)$, and its generating function is given by
\begin{equation}\label{gf1}
  \sum_{n=0}^\infty \overline{p}(n)\,q^n 
  \;=\; \frac{(-q;q)_\infty}{(q;q)_\infty},
\end{equation}
where $(x;q)_\infty = \displaystyle
\prod_{k=0}^{\infty}(1 - x\,q^k)$ for $|q|<1$.

An overpartition pair of $n$ is a pair of overpartitions $(\lambda_1, \lambda_2)$, where the sum of all parts is $n$. Let $\overline{pp}(n)$ 
denote the number of such pairs, whose generating function given by 
\begin{equation}\label{E29}  
\sum_{n=0}^{\infty} \overline{pp}(n) q^n =\frac{(-q;q)^2_\infty}{(q;q)^2_\infty}.
\end{equation}

The arithmetic properties of $\overline{p}(n)$ and $\overline{pp}(n)$ have been the subject of numerous recent studies \cite{BrL, ChL, HS, K, M, T}.

Recently, Lin \cite{BL} and Adiga and Ranganatha \cite{AR} examined the arithmetic properties of the number of overpartition pairs of $n$ restricted to odd parts, $\overline{pp}_{o}(n)$. For instance, authors \cite[Theorem 1.2]{AR} have proved that 
\begin{equation}
	\overline{pp}_{o}(2^{\beta+2} n)\equiv 0 \pmod{2^{3\beta +5}}\label{eq6}
\end{equation}
for all integers $n\geq 1$ and $\beta \geq 3$.

Most recently, Naika and Shivashankar \cite{MS} established infinite families of congruences modulo $3,8,16,36,48,$ and $96$ for $\overline{B}_{3}(n)$ and modulo $3,16,64,$ and $96$ for $\overline{B}_{4}(n)$,
where $\overline{B}_{k}(n)$ counts the number of $k$--regular overpartition pairs of $n$. The corresponding generating function is given by
\begin{equation*}  
\sum_{n=0}^{\infty} \overline{B}_{k}(n) q^n =\frac{(-q;q)^2_\infty(q^k;q^k)^2_\infty}{(q;q)^2_\infty (-q^k;q^k)^2_\infty}.
\end{equation*}
Rahman and Saikia \cite{RN} further investigated congruences modulo powers of 2 for $\overline{B}_{3\alpha}(n)$, $\overline{B}_{4\alpha}(n)$ and $\overline{B}_{5\alpha}(n)$, where $\alpha$ is a positive integer. Shivashankar and Gireesh \cite{SG} established several infinite families of congruences for $\overline{B}_{k}(n)$ when $k = 3, 4, 5, 8$.

In this paper, we extend these investigations by proving infinite families of congruences modulo powers of 2 for $2^\alpha$--regular overpartition pairs, $\overline{B}_{2^\alpha}(n)$.
The main results of the paper are as follows.

\begin{theorem}\label{T1}
For any integers $n\geq0, \alpha\geq1$ and $\beta \geq 0$, we have
\begin{equation}\label{eq1}
\overline{B}_{2^{\alpha}}(2^{\alpha}(n+1)) \equiv 0\pmod{2^3}
\end{equation}
and
\begin{equation}\label{eq2}
\overline{B}_{2^{\alpha}}(2^{\alpha+\beta+1}(n+1)) \equiv 0\pmod{2^{3\beta+5}}.
\end{equation}
\end{theorem}
\noindent Note that \eqref{eq6} is a special case of \eqref{eq2} with $\alpha=1$.

For any integer $n$, let $\sigma(n)$ denote the sum of all positive divisors of $n$ and let $\nu_p(n)$ be the highest power of the prime $p$ that divides $n$.
\begin{theorem}\label{T2}
For any $n\geq 0, \, \alpha, \beta \geq1$, we have
\begin{align}
\overline{B}_{2^{\alpha}}(2^{\alpha}(2n+1)) &\equiv y_{2, 1}^{(\alpha)} \sigma(2n+1) \pmod{2^6}, \label{eq3} \\
\overline{B}_{2^{\alpha}}(2^{\alpha+1}(n+1)) &\equiv y_{1,1}^{(\alpha)} (n+1)^3 \sigma(n+1) \pmod{2^9} \label{eq4}
\end{align}
and
\begin{equation}\label{eq5}
\overline{B}_{2^{\alpha}}(2^{\alpha+\beta+1}(n+1)) \equiv z_{\beta,1}^{(\alpha)} (n+1)^3 \sigma(n+1) \pmod{2^{3\beta+9}}.
\end{equation}
\end{theorem}
\begin{coro}\label{Co1}
Suppose that $n\geq0$, $1\leq j <4$, and $p$ is a prime. If $p \equiv -1 \pmod{2^j}$ and $\nu_p(2n+1) \equiv -1 \pmod{2^{4-j}}$, then for any integers $\alpha \geq 1$,
\begin{equation}
\overline{B}_{2^{\alpha}}(2^{\alpha}(2n+1)) \equiv 0 \pmod{2^6}. \label{EE17} 
\end{equation}
\end{coro}
\begin{coro}\label{Co2}
	Let $n\geq 0$, $1\leq j <5$, and $p$ be a prime such that $p \equiv -1 \pmod{2^j}$ and $\nu_p(n+1) \equiv -1 \pmod{2^{5-j}}$. Then for any integers $\alpha\geq1$ and $\beta \geq0$, 
	\begin{equation}
		\overline{B}_{2^{\alpha}}(2^{\alpha+\beta+1}(n+1)) \equiv 0 \pmod{2^{3\beta+9}}. \label{EE19}
	\end{equation}
\end{coro}
\begin{coro}\label{Co3}
	If $n\geq0$ and $p$ is a  prime such that $p \nmid (2n+1)$ and $p\equiv -1 \pmod{2^{4-j}}$, $1\leq j < 4$, then, for any integers $\alpha, k\geq1$, 
	\begin{equation}
		\overline{B}_{2^{\alpha}}(2^{\alpha}p^{2^j k-1}(2n+1))\equiv 0 \pmod{2^6}. \label{EE20}
	\end{equation}
\end{coro}

\begin{coro}\label{Co4}
If $n\geq0$ and $p$ is a prime such that $p \nmid (n+1)$ and $p\equiv -1 \pmod{2^{5-j}}$, $1\leq j < 5$, then, for any integers $\alpha, k\geq1$, and $\beta \geq 0$,
\begin{equation}
\overline{B}_{2^{\alpha}}(2^{\alpha+\beta+1}p^{2^j k-1}(n+1))\equiv 0 \pmod{2^{3\beta+9}}. \label{EE22}
\end{equation}
\end{coro}

\section{Preliminaries} \label{S2}
In this section, we present some preliminary results, which are useful in proving our main results.

\noindent Ramanujan's theta function $f(a,b)$, defined by 
\begin{equation*} \label{E1}
	f(a,b) = \sum\limits_{n=-\infty}^{\infty} a^{n(n+1)/2}b^{n(n-1)/2}, \qquad |ab|<1.
\end{equation*}\label{E2}
Jacobi's triple product identity takes the form
\begin{equation*}\label{E3}
	f(a,b) = (-a;ab)_{\infty} (-b;ab)_{\infty} (ab;ab)_{\infty}.
\end{equation*}
Some special cases of $f(a,b)$ are given by 
\begin{equation*}\label{E4}
\psi(q):= f(q,q^3)=\sum\limits_{n=0}^{\infty} q^{n(n+1)/2} = \dfrac{(q^2;q^2)_\infty}{(q;q^2)_\infty}
\end{equation*} and

\begin{equation*}\label{E5}
	\varphi(q):= f(q,q)=\sum\limits_{n=0}^{\infty} q^{n^2} = (-q;q^2)^2_\infty(q^2;q^2)_\infty.
\end{equation*}
Consider the identities \cite[Entry 3, p. 40]{BCB3} by Ramanujan,
\begin{lemma} \label{L1}
	\begin{align}
		\varphi(-q^2)^2&=\varphi(q)\varphi(-q), \label{E6}\\
		\psi(q)^2&=\varphi(q)\psi(q^2), \label{E7}\\
		\varphi(q)&=\varphi(q^4)+2q\psi(q^8), \label{E8}\\
		\varphi(q)^2&=\varphi(q^2)^2+4q\psi(q^4)^2. \label{E9}
	\end{align}
\end{lemma}
We have the Huffing operator $H$ defined by:
\begin{equation*}
H\lb (\sum\limits_{n=0}^{\infty} a_nq^n\rb )=\sum\limits_{n=0}^{\infty} a_{2n}q^{2n}. 
\end{equation*}
\begin{lemma}\cite[(2.20)]{HB}\label{Le22}
	For any $j\geq1$, we have
\begin{equation} \label{P6}
	H\left(\dfrac{\varphi(q^2)^{2j}}{\varphi(-q)^{2j}}\right) = \sum_{k=1}^j m(j, k) \dfrac{\varphi(q^2)^{4k}}{\varphi(-q^2)^{4k}},
\end{equation}
where the matrix $M=\{m(j,k)\}_{j,k\ge1}$ is defined as follows: 
\begin{equation*}
m(1,1)=1, \, m(2,1)=-1,\, m(2,2)=2,
\end{equation*}
for all $k \geq 2$
\begin{equation*}
m(1, k)=0, \, m(2, k+1)=0,
\end{equation*}
and for all $k \geq 2$ and $j \geq 3$
\begin{equation*}
m(j, 1)=0, \, m(j, k)=2\,m(j-1, k-1)-m(j-2, k-1).
\end{equation*}
\end{lemma}
The first five rows of $M$ are given by
\begin{equation*}
	\begin{pmatrix} 
		1&0&0&0&0&0&\cdots\\
		-1&2&0&0&0&0&\cdots\\
		0&-3&2^2&0&0&0&\cdots\\
		0&1&-2^3&2^3&0&0&\cdots\\
		0&0&5&-2^2\times5&2^4&0&\cdots
	\end{pmatrix}.
\end{equation*}
Using the definition of $M$ and induction, it is easy to prove that 
\begin{enumerate}
\item \label{MR1} $m(j, k)=0$ for all $k>j$.
\item \label{MR2} $m(j, k)=0$ for all $j>2k$.
\item \label{MR3} $\displaystyle \sum_{k=1}^\infty m(j, k) = 1$ for all $j\geq1$.
\item \label{MR4} $m(j, j)=2^{j-1}$ for all $j\geq 1$.
\end{enumerate}

\begin{lemma}\cite[(2.31)]{HB} \label{Le24}
For any $j\geq1$, we have
\begin{equation} \label{P12}
	H\left(q^j\dfrac{\psi(q^4)^{2j}}{\varphi(-q)^{2j}}\right) = \sum_{k=1}^j r(j, k)  q^{2k}\dfrac{\psi(q^4)^{4k}}{\varphi(-q^2)^{4k}},
\end{equation}
where the matrix $R=\{r(j, k)\}_{j,\,k\ge1}$ is defined as follows:
\begin{equation*}
	r(1,1)=4, \, r(2,1)=1,\, r(2,2)=32,
\end{equation*}
for all $k \geq 2$
\begin{equation*}
	r(1, k)=0, \, r(2, k+1)=0,
\end{equation*}
and for all $k \geq 2$ and $j \geq 3$
\begin{equation*}
	r(j, 1)=0, \, r(j, k)=8\,r(j-1, k-1)+r(j-2, k-1).
\end{equation*}
\end{lemma}
The first five rows of $R$ are given by
\begin{equation*}
	\begin{pmatrix} 2^2&0&0&0&0&\cdots\\
		1&2^5&0&0&0&\cdots\\
		0&2^2\times 3&2^{8}&0&0&\cdots\\
		0&1&2^7&2^{11}&0&\cdots\\
		0&0&2^2\times 5&2^{8}\times 5&2^{14}&\cdots
	\end{pmatrix}.
\end{equation*}
\begin{lemma}\cite[(2.10)]{HRB} \label{Le25}
For any $j\geq1$, we have
\begin{equation} \label{P8}
H\left(q^j\dfrac{\psi(q^2)^{4j}}{\varphi(-q)^{4j}}\right) = \sum_{k=1}^j n(j, k)  q^{2k}\dfrac{\psi(q^2)^{8k}}{\varphi(-q^2)^{8k}},
\end{equation}
where the matrix $N=\{n(j, k)\}_{j,\,k\ge1}$ is defined as follows:
\begin{equation*}
	n(1,1)=8, \, n(2,1)=1,\, n(2,2)=128,
\end{equation*}
for all $k \geq 2$
\begin{equation*}
	n(1, k)=0, \, n(2, k+1)=0,
\end{equation*}
and for all $k \geq 2$ and $j \geq 3$
\begin{equation*}
	n(j, 1)=0, \, n(j, k)=16\,n(j-1, k-1)+n(j-2, k-1).
\end{equation*}
\end{lemma}
The first five rows of $N$ are given by
\begin{equation*}
	\begin{pmatrix} 2^3&0&0&0&0&\cdots\\
		1&2^7&0&0&0&\cdots\\
		0&2^3\times 3&2^{11}&0&0&\cdots\\
		0&1&2^9&2^{15}&0&\cdots\\
		0&0&2^3\times 5&2^{11}\times 5&2^{19}&\cdots
	\end{pmatrix}.
\end{equation*}
Using the definition of $N$ and induction, it is easy to prove that 
\begin{enumerate}
	\setcounter{enumi}{4}
	\item \label{NR1} $n(j, k)=0$ for all $k>j$.
	\item \label{NR2} $n(j, k)=0$ for all $j>2k$.
	\item \label{NR3} $n(2j, j)=1$ for all $j\geq 1$.
	\item \label{NR4} $n(j, j)=2^{4j-1}$ for all $j\geq 1$.
\end{enumerate}

\begin{lemma} 
For any $j\geq 1$, we have
	\begin{equation}\label{E11}
		H\lb ( \dfrac{\varphi(q)^{4j-2}}{\varphi(-q)^{4j}}\rb)=\sum\limits_{k=1}^{2j} m(4j-1, k+2j-1)\,\, \dfrac{\varphi(q^2)^{4k-2}}{\varphi(-q^2)^{4k}}.
	\end{equation}
\end{lemma}
\begin{proof}
By replacing $j$ by $4j-1$ in \eqref{P6} and using the fact that $m(4j-1, k)=0$ for $1\leq k <2j$, we obtain
\begin{equation*} 
	H\left(\dfrac{\varphi(q^{2})^{8j-2}}{\varphi(-q)^{8j-2}}\right) = \sum_{k=1}^{2j} m(4j-1, k+2j-1) \dfrac{\varphi(q^2)^{4k+8j-4}}{\varphi(-q^2)^{4k+8j-4}},
\end{equation*}
By rewriting the above equation and using \eqref{E6}, we obtain \eqref{E11}.
\end{proof}

\begin{lemma} \cite[Eqns (3.23) and (3.24)]{HB1} 
For any $j\geq 1$, we have
\begin{equation}\label{E16}
		H\lb ( \dfrac{\varphi(q)^{4j-2}}{\varphi(-q)^{4j-2}}\rb)=\sum_{k=1}^{2j} r(4j-2, k+2j-2) \,q^{2k-2}\dfrac{\psi(q^4)^{4k-4}}{\varphi(-q^2)^{4k-4}}
\end{equation}
and
	\begin{equation}\label{E17}
		H\lb (q\, \dfrac{\varphi(q)^{4j-2}}{\varphi(-q)^{4j-2}}\rb)=\sum\limits_{k=1}^{j} n(2j-1, k+j-1)\,q^{2k} \dfrac{\psi(q^2)^{8k-4}}{\varphi(-q^2)^{8k-4}}.
	\end{equation}
\end{lemma}

\begin{lemma} \cite[Eqn. (3.30)]{HB1}
For any integer $j\geq 1$, we have
\begin{equation}\label{E12}
H\lb (q^{j}\frac{\psi(q)^{8j}}{\varphi(-q)^{8j}}\rb)=\sum_{k=1}^{2j} n(3j, k+j)\,q^{2k}\, \frac{\psi(q^2)^{8k}}{\varphi(-q^2)^{8k}}.
\end{equation}
\end{lemma}

\section{Generating Functions} \label{S3}
In this section, we establish generating functions for $ \overline{pp}(n)$ and  $ \overline{B}_{2^\alpha}(n)$, in arithmetic progressions.
\begin{lemma}\label{L4}
For any integer $\alpha>0$,
\begin{equation}\label{E26}
\sum\limits_{n=0}^{\infty} \overline{pp}(2^{\alpha}n) q^{n} = \sum\limits_{k=1}^{2^{\alpha-1}} x_{\alpha, k} \dfrac{\varphi(q)^{4k-2}}{\varphi(-q)^{4k}}.
\end{equation}
The coefficient vector 
$\textbf{x}_\alpha = (x_{\alpha, 1},\ x_{\alpha, 2},\ x_{\alpha, 3}, \ldots)$ are defined by 
\begin{equation}\label{E261}
\textbf{x}_1 = (1,\ 0,\ 0, \ldots )
\end{equation}
and for all $\alpha \geq 1$, 
\begin{equation}\label{E27}
\textbf{x}_{\alpha+1} = \textbf{x}_{\alpha} A,
\end{equation}
where $A=a(j,k)_{j,k\geq 1}$ is a submatrix of $M$,
\begin{equation*}\label{E28}
	a(j,k) = m(4j-1, k+2j-1).
\end{equation*}
The first three cases are given by
\begin{align*}
	\bf{x_1} &= (1,\ 0,\ 0,\ 0, \ldots),\\
	\bf{x_2} &= (-3,\ 4,\ 0,\ 0, \ldots),\\
	\bf{x_3} &= (-19,\ 2^2\times 53,\ -2^6\times 7,\ 2^8,\ 0,\ 0 \ldots).
\end{align*}
\end{lemma}
\begin{proof}
From \eqref{E29}, we have
\begin{equation*}
\sum\limits_{n=0}^{\infty} \overline{pp}(2n) q^n =\frac{\varphi(q)^2}{\varphi(-q^2)^4}.
\end{equation*}
	Using \eqref{E9} and extracting the terms that involve even powers of $q$, we obtain the case $\alpha=1$ of \eqref{E26}. 
    
Now, suppose \eqref{E26} holds for some $\alpha \geq 1$. Then
\begin{equation*}
\sum\limits_{n=0}^{\infty} \overline{pp}(2^{\alpha}n) q^{n} = \sum\limits_{k=1}^{2^{\alpha-1}} x_{\alpha, k} \dfrac{\varphi(q)^{4k-2}}{\varphi(-q)^{4k}}.
\end{equation*}
    If we apply the operator $H$ to the above equation and use \eqref{E11}, we obtain
	\begin{align*}\label{E30}
		\sum\limits_{n=1}^{\infty} \overline{pp}(2^{\alpha+1}n) q^{2n} &= \sum\limits_{k=1}^{2^{\alpha-1}} x_{\alpha, k} H \lb( \dfrac{\varphi(q)^{4k-2}}{\varphi(-q)^{4k}}\rb)\\
		&= \sum\limits_{k=1}^{2^{\alpha-1}} x_{\alpha, k} 
		\sum\limits_{j=1}^{2k} m(4k-1, j+2k-1)\,\, \dfrac{\varphi(q^2)^{4j-2}}{\varphi(-q^2)^{4j}}\\
		&=  \sum\limits_{j=1}^{2^{\alpha}}  
		\lb( \sum\limits_{k=1}^{2^{\alpha -1}}x_{\alpha, k} \, m(4k-1, j+2k-1)\rb) \dfrac{\varphi(q^2)^{4j-2}}{\varphi(-q^2)^{4j}}\\
		&= \sum\limits_{j=1}^{2^{\alpha}}x_{\alpha+1, j} \,\, \dfrac{\varphi(q^2)^{4j-2}}{\varphi(-q^2)^{4j}}
	\end{align*}
	which is \eqref{E26} with $\alpha+1$ for $\alpha$.
\end{proof}
\begin{lemma}\label{L5}
	For any integer $\alpha\geq 1$,
	\begin{equation}\label{E31}
		\sum\limits_{n=0}^{\infty} \overline{B}_{2^{\alpha}}(2^{\alpha}n) q^{n} = \sum\limits_{k=1}^{2^{\alpha-1}} x_{\alpha, k} \dfrac{\varphi(q)^{4k-2}}{\varphi(-q)^{4k-2}}.
	\end{equation}
\end{lemma}
\begin{proof}
	It is evident that 
	\begin{align*}\label{E32}
		\sum\limits_{n=0}^{\infty} \overline{B}_{2^{\alpha}}(n) q^{n} &= \dfrac{\varphi(-q^{2^\alpha})^{2}}{\varphi(-q)^2}\\
		&=\varphi(-q^{2^\alpha})^{2} \sum\limits_{m=0}^{\infty} \overline{pp}(m) q^{m}.
	\end{align*}
	Extracting the terms involving $q^{2^{\alpha}}$ on both sides of the above equation, we obtain
	\begin{equation*}\label{E33}
		\sum\limits_{n=0}^{\infty} \overline{B}_{2^{\alpha}}(2^{\alpha}n) q^{2^{\alpha}n} = \varphi(-q^{2^\alpha})^{2} \sum\limits_{m=0}^{\infty} \overline{pp}(2^{\alpha}m) q^{2^{\alpha}m}.
	\end{equation*}
	On replacing $q^{2^{\alpha}}$ by $q$ and using \eqref{E26}, we deduce that 
	\begin{equation*}\label{E34}
	\sum\limits_{n=0}^{\infty} \overline{B}_{2^{\alpha}}(2^{\alpha}n) q^{n} = \varphi(-q)^{2} \sum\limits_{k=1}^{2^{\alpha-1}} x_{\alpha, k} \dfrac{\varphi(q)^{4k-2}}{\varphi(-q)^{4k}},
\end{equation*}
which is \eqref{E31}.
\end{proof}
\begin{lemma}\label{L6}
	For any integer $\alpha\geq 1$, we have 
	\begin{equation}\label{E35}
		\sum\limits_{k=1}^{2^{\alpha -1 }} x_{\alpha, k} =1.
	\end{equation}
\end{lemma}
\begin{proof}
	The result is true for $\alpha =1$. Suppose that it holds for any $\alpha \geq 1$. Then
	\begin{align*}\label{E36}
		\sum_{k=1}^{2^{\alpha -1 }} x_{\alpha+1, k} &= \sum\limits_{k=1}^{2^{\alpha -1 }} \lb( \sum\limits_{j=1}^{2^{\alpha -1 }} x_{\alpha, j} \, m(4j-1, k+2j-1)\rb)\\
		&= \sum\limits_{j=1}^{2^{\alpha -1 }} x_{\alpha, j} \lb( \sum\limits_{k=1}^{2^{\alpha -1 }}  m(4j-1, k+2j-1)\rb)\\
		&= \sum\limits_{j=1}^{2^{\alpha -1 }} x_{\alpha, j} =1.
	\end{align*}
	We have used \eqref{E27} in first equality and the property \eqref{MR3} of matrix $M$ in the second equality. This completes the proof by induction.
\end{proof}
\begin{lemma}\label{L7}
	For any integer $\alpha>0$,
	\begin{equation}\label{E37}
		\sum_{n=0}^{\infty} \overline{B}_{2^{\alpha}}(2^{\alpha+1}(n+1)) q^{n} =\sum\limits_{k=1}^{2^{\alpha}-1} y_{1, k}^{(\alpha)} \, q^{k-1}  \dfrac{\psi(q^2)^{4k}}{\varphi(-q)^{4k}}
	\end{equation}
	and
		\begin{equation}\label{E38}
		\sum\limits_{n=0}^{\infty} \overline{B}_{2^{\alpha}}(2^{\alpha}(2n+1)) q^{n} = \sum\limits_{k=1}^{2^{\alpha-1}} y_{2, k}^{(\alpha)} \, q^{k-1}  \dfrac{\psi(q)^{8k-4}}{\varphi(-q)^{8k-4}}.
	\end{equation}
The coefficient vectors are defined as follow:
\begin{equation*}
y_{1, k}^{(\alpha)} = \sum_{j=1}^{2^{\alpha-1}} x_{\alpha, j} \, r(4j-2, k+2j-1)
\end{equation*}
and
\begin{equation*}
y_{2, k}^{(\alpha)} = \sum_{j=1}^{2^{\alpha-1}} x_{\alpha, j}\, n(2j-1, k+j-1).
\end{equation*}
\end{lemma}
\begin{proof}
If we apply $H$ to \eqref{E31} and use \eqref{E16}, we obtain
\begin{align*}
\sum_{n=0}^{\infty} \overline{B}_{2^{\alpha}}(2^{\alpha+1}n) q^{2n} &= \sum_{k=1}^{2^{\alpha-1}} x_{\alpha, k} \,H\lb(\frac{\varphi(q)^{4k-2}}{\varphi(-q)^{4k-2}}\rb)\\&
=\sum_{k=1}^{2^{\alpha-1}} x_{\alpha, k} \sum_{j=1}^{2k} r(4k-2, j+2k-2)\,q^{2j-2} \frac{\psi(q^4)^{4j-4}}{\varphi(-q^2)^{4j-4}}\\&
=\sum_{j=1}^{2^\alpha} \lb(\sum_{k=1}^{2^{\alpha-1}} x_{\alpha, k} \,r(4k-2, j+2k-2)\rb)\,q^{2j-2}\frac{\psi(q^4)^{4j-4}}{\varphi(-q^2)^{4j-4}}.
\end{align*}
From Lemma \ref{L6} and the fact that $r(2, 1)=1$, we can rewrite the above equation as
\begin{align*}
\sum_{n=0}^{\infty} \overline{B}_{2^{\alpha}}(2^{\alpha+1}n) q^{2n} &= 1+\sum_{j=1}^{2^\alpha-1} \lb(\sum_{k=1}^{2^{\alpha-1}} x_{\alpha, k} \,r(4k-2, j+2k-1)\rb)\,q^{2j}\frac{\psi(q^4)^{4j}}{\varphi(-q^2)^{4j}} \\&
=1+\sum_{j=1}^{2^\alpha-1} y_{1, j}^{(\alpha)}\,q^{2j}\frac{\psi(q^4)^{4j}}{\varphi(-q^2)^{4j}},
\end{align*}
which yields \eqref{E37}.

From \eqref{E31}, we have
\begin{equation*}
		\sum\limits_{n=0}^{\infty} \overline{B}_{2^{\alpha}}(2^{\alpha}n) q^{n+1} = \sum\limits_{k=1}^{2^{\alpha-1}} x_{\alpha, k} \,q\frac{\varphi(q)^{4k-2}}{\varphi(-q)^{4k-2}}
	\end{equation*}
    If we now employ operator $H$ on both sides of the above equation and use \eqref{E17}, we obtain
\begin{align*}
	\sum_{n=0}^{\infty} \overline{B}_{2^{\alpha}}(2^{\alpha}(2n+1)) q^{2n+2} &= \sum\limits_{k=1}^{2^{\alpha-1}} x_{\alpha, k} \,H\lb(q\frac{\varphi(q)^{4k-2}}{\varphi(-q)^{4k-2}}\rb)\\&
    =\sum_{k=1}^{2^{\alpha-1}} x_{\alpha, k} 
    \sum_{j=1}^{k} n(2k-1, j+k-1)\,q^{2j} \dfrac{\psi(q^2)^{8j-4}}{\varphi(-q^2)^{8j-4}} \\&
    =\sum_{j=1}^{2^{\alpha-1}} \lb(\sum_{k=1}^{2^{\alpha-1}} x_{\alpha, k}\, n(2k-1, j+k-1)\rb)q^{2j} \dfrac{\psi(q^2)^{8j-4}}{\varphi(-q^2)^{8j-4}}\\&
    =\sum_{j=1}^{2^{\alpha-1}} y_{2, j}^{(\alpha)}
    \, q^{2j} \dfrac{\psi(q^2)^{8j-4}}{\varphi(-q^2)^{8j-4}}.
\end{align*}
If we replace $q^2$ in place of $q$, we obtain \eqref{E38}.
\end{proof}
\begin{lemma}\label{L8}
For any integers $\alpha, \beta \geq1$,
\begin{equation}\label{E32}
\sum_{n=0}^\infty \overline{B}_{2^{\alpha}}(2^{\alpha+\beta+1}(n+1)) q^n = \sum_{k=1}^{2^{\beta-1}(2^\alpha-1)} z_{\beta, k}^{(\alpha)} \, q^{k-1} \frac{\psi(q)^{8k}}{\varphi(-q)^{8k}},
\end{equation}
where the coefficient vectors are defined as follow:
\begin{equation*}
z_{1, j}^{(\alpha)} = \sum_{k=1}^{2^\alpha-1} y_{1, k}^{(\alpha)}\,\, n(k, j)
\end{equation*}
and
\begin{equation*}
z_{\beta+1, j}^{(\alpha)} = \sum_{k=1}^{2^{\beta-1}(2^\alpha-1)} z_{\beta, k}^{(\alpha)}\,\, n(3k, k+j).
\end{equation*}
\end{lemma}
\begin{proof}
From \eqref{E37}, we have
\begin{equation*}
\sum_{n=0}^{\infty} \overline{B}_{2^{\alpha}}(2^{\alpha+1}(n+1)) q^{n+1} = \sum_{k=1}^{2^{\alpha}-1} y_{1, k}^{(\alpha)} \, q^{k}  \frac{\psi(q^2)^{4k}}{\varphi(-q)^{4k}}.
\end{equation*}
If we now apply $H$ and use \eqref{P8}, we obtain
\begin{align*}
\sum_{n=0}^{\infty} \overline{B}_{2^{\alpha}}(2^{\alpha+2}(n+1)) q^{2n+2} & 
=\sum_{k=1}^{2^{\alpha}-1} y_{1, k}^{(\alpha)} \sum_{j=1}^k n(k, j) q^{2j} \frac{\psi(q^2)^{8j}}{\varphi(-q^2)^{8j}} \\&
=\sum_{j=1}^{2^{\alpha}-1} \lb( \sum_{k=1}^{2^{\alpha}-1} y_{1, k}^{(\alpha)}\, n(k, j)\rb) q^{2j} \frac{\psi(q^2)^{8j}}{\varphi(-q^2)^{8j}}\\&
=\sum_{j=1}^{2^{\alpha}-1} z_{1, j}^{(\alpha)} q^{2j} \frac{\psi(q^2)^{8j}}{\varphi(-q^2)^{8j}}
\end{align*}
which yields
\begin{equation*}
\sum_{n=0}^{\infty} \overline{B}_{2^{\alpha}}(2^{\alpha+2}(n+1)) q^{n+1} = \sum_{j=1}^{2^{\alpha}-1} z_{1, j}^{(\alpha)} q^{j} \frac{\psi(q)^{8j}}{\varphi(-q)^{8j}}
\end{equation*}
This implies that \eqref{E27} is true for $\beta = 1$. Suppose that \eqref{E27} is true for an integer $\beta \geq 1$. Then
\begin{equation*}
\sum_{n=0}^\infty \overline{B}_{2^{\alpha}}(2^{\alpha+\beta+1}(n+1)) q^{n+1} = \sum_{k=1}^{2^{\beta-1}(2^\alpha-1)} z_{\beta, k}^{(\alpha)} \, q^{k} \frac{\psi(q)^{8k}}{\varphi(-q)^{8k}}
\end{equation*}
If we apply $H$ on both sides and use \eqref{E12}, we obtain
\begin{align*}
\sum_{n=0}^\infty \overline{B}_{2^{\alpha}}(2^{\alpha+\beta+2}(n+1)) q^{2n+2} &
= \sum_{k=1}^{2^{\beta-1}(2^\alpha-1)} z_{\beta, k}^{(\alpha)}\,  \sum_{j=1}^{2k} n(3k, j+k)\,q^{2j}\, \frac{\psi(q^2)^{8j}}{\varphi(-q^2)^{8j}}\\&
=\sum_{j=1}^{2^{\beta}(2^\alpha-1)}\, \sum_{k=1}^{2^{\beta-1}(2^\alpha-1)} z_{\beta, k}^{(\alpha)} \,n(3k, j+k) \,q^{2j}\,\frac{\psi(q^2)^{8j}}{\varphi(-q^2)^{8j}} \\&
=\sum_{j=1}^{2^{\beta}(2^\alpha-1)} z_{\beta+1, j}^{(\alpha)} \,q^{2j}\,\frac{\psi(q^2)^{8j}}{\varphi(-q^2)^{8j}},
\end{align*}
which is \eqref{E32} with $\beta+1$ in place of $\beta$. This completes the proof.
\end{proof}
\section{Congruences}
For a positive integer $n$, let $\nu_2(n)$ be the highest power of 2 that divides $n$, and define $\nu_2(0)=+\infty$.
In this section, we prove some lemmas which play a crucial role in proving our main results.

\begin{lemma}
For any integers $j, k\geq 1$, we have
\begin{align}
\nu_2(m(j, k)) &\geq 2k-j-1, \label{L31}\\
\nu_2(r(j, k)) &\geq 6k-3j-1, \label{L32}
\end{align}
and 
\begin{equation}
\nu_2(n(j, k)) \geq 8k-4j-1. \label{L33}
\end{equation}
\end{lemma}
\begin{proof}
The proof follows from the definition of matrices $M, R, N$ and induction.
\end{proof}
\begin{lemma}
For any integers $\alpha, k\geq1$, we have
\begin{equation}\label{L34}
\nu_2(x_{\alpha, k}) \geq 2k-2.
\end{equation}
\end{lemma}
\begin{proof}
The proof follows from \eqref{E261}, \eqref{E27} and \eqref{L31}.
\end{proof}
\begin{lemma} \label{L9}
For any integers $\alpha, k\geq1$, we have
\begin{equation}\label{L35}
\nu_2\left(y_{1, k}^{(\alpha)}\right) \geq 6k-1
\end{equation}
and
\begin{equation}\label{L36}
\nu_2\left(y_{2, k}^{(\alpha)}\right) \geq 8k-5.
\end{equation}
\end{lemma}
\begin{proof}
For Lemma \ref{L7}, we have
\begin{align*}
\nu_2\left(y_{1, k}^{(\alpha)}\right) &\geq \min_{j\geq 1} (\nu_2(x_{\alpha, j})+\nu_2 (r(4j-2, k+2j-1))) \\& \geq 6k-1
\end{align*}
and
\begin{align*}
\nu_2\left(y_{2, k}^{(\alpha)}\right) &\geq \min_{j\geq 1} (\nu_2(x_{\alpha, j})+\nu_2 (n(2j-1, k+j-1))) \\& \geq 8k-5.
\end{align*}

\end{proof}
\begin{lemma}
For any integer $\alpha, \beta, k\geq1$, we have
\begin{equation}\label{L37}
\nu_2(z_{\beta, k}^{(\alpha)}) \geq 8k+3\beta-3.
\end{equation}
\end{lemma}
\begin{proof}
From Lemmas \ref{L8} and \ref{L9} we have
\begin{align*}
\nu_2\left(z_{1, k}^{(\alpha)}\right) &\geq \min_{j\geq1} (\nu_2(y_{1,j}^{(\alpha)}) + \nu_2(n(j, k)))\\
& \geq 8k.
\end{align*}
Now suppose \eqref{L37} is true for some $\beta \geq1$. Then by \eqref{L33},
\begin{align*}
\nu_2\left(z_{\beta+1, k}^{(\alpha)}\right) & \geq \min_{j\geq1} (\nu_2(z_{\beta,j}^{(\alpha)}) + \nu_2(n(3j, k+j))) \\
& \geq 8k+3\beta,
\end{align*}
which is \eqref{L37} with $\beta+1$ for $\beta$.
\end{proof}
We are now ready to prove Theorems \ref{T1}, \ref{T2} and corollaries \ref{Co1} - \ref{Co4}.
Let $t_k(n)$ denote the number of representations of a nonnegative integer $n$ as a sum of $k$ triangular numbers.  Then
\begin{equation}\label{E40}
	\sum_{n=0}^\infty t_k(n) q^n = \psi(q)^k.
\end{equation}
\subsection*{Proof of Theorem \ref{T1}} By the binomial theorem, we have
\[\frac{\varphi(q)^{2j}}{\varphi(-q)^{2j}} \equiv 1 \pmod{2^{3}}.\]
From \eqref{E31} and \eqref{E35}, we see that 
\begin{equation}\label{eeqn1}
\sum_{n=0}^\infty \overline{B}_{2^{\alpha}}(2^{\alpha}n) \,q^{n} \equiv1 \pmod{2^3}.
\end{equation}
Congruence \eqref{eq1} follows from \eqref{eeqn1}. 

In view of \eqref{E37} and \eqref{L35}, we have
\begin{equation*}
\sum_{n=0}^{\infty} \overline{B}_{2^{\alpha}}(2^{\alpha+1}(n+1)) q^{n} \equiv y_{1, k}^{(\alpha)}  \frac{\psi(q^2)^{4}}{\phi(-q)^{4}} \pmod{2^{11}}.
\end{equation*}
By the binomial theorem, it is easy that
\[\frac{\psi(q^2)^{4}}{\phi(-q)^{4}} \equiv \psi(q)^8 \pmod{2^4}.\]
Thus,
\begin{equation}\label{eeq21}
\sum_{n=0}^{\infty} \overline{B}_{2^{\alpha}}(2^{\alpha+1}(n+1)) q^{n} \equiv y_{1, k}^{(\alpha)}  \psi(q)^{8} \pmod{2^{9}}.
\end{equation}
In view of \eqref{E32} and \eqref{L37}, 
\begin{equation*}
\sum_{n=0}^\infty \overline{B}_{2^{\alpha}}(2^{\alpha+\beta+1}(n+1)) q^n \equiv  z_{\beta, 1}^{(\alpha)} \frac{\psi(q)^{8}}{\varphi(-q)^{8}}  \pmod{2^{3\beta+13}}.
\end{equation*}
Since
\[\frac{1}{\varphi(-q)^8}\equiv 1 \pmod{2^4},\]
we have
\begin{equation}\label{eeq22}
\sum_{n=0}^\infty \overline{B}_{2^{\alpha}}(2^{\alpha+\beta+1}(n+1)) q^n \equiv  z_{\beta, 1}^{(\alpha)} \psi(q)^{8}  \pmod{2^{3\beta+9}}.
\end{equation}
Congruence \eqref{eq2} follows from \eqref{eeq21}, \eqref{eeq22} and the fact that $\nu_2\lb(y_{1, 1}^{(\alpha)}\rb)=5$ and $\nu_2\lb(z_{\beta, 1}^{(\alpha)}\rb)=3\beta+5$ for each $\beta \geq 1$.
\subsection*{Proof of Theorem \ref{T2}}
From \eqref{E38}, \eqref{L36}, and 
\[\frac{1}{\varphi(-q)^4}\equiv 1 \pmod{2^3},\]
we have
\begin{equation}\label{E39}
		\sum\limits_{n=0}^{\infty} \overline{B}_{2^{\alpha}}(2^{\alpha}(2n+1)) q^{n} \equiv y_{2, 1}^{(\alpha)} \psi(q)^4 \pmod{2^6}.
	\end{equation}
    
From \cite[Theorem 3.6.3]{BC}, we have
\begin{equation}\label{E41}
t_4(n) = \sigma(2n+1).
\end{equation}
Congruence \eqref{eq3} follows from \eqref{E40}, \eqref{E39} and \eqref{E41}.
From \cite[Equation 4.32]{HB1}, we have
\begin{equation}\label{E42}
t_8(n) = (n+1)^3 \sigma(n+1) \pmod{2^4}. 
\end{equation}
Similarly, congruences \eqref{eq4} and \eqref{eq5} follow from \eqref{eeq21}, \eqref{eeq22}, \eqref{L35}, \eqref{L37}, \eqref{E40} and \eqref{E42}.

\subsection*{Proofs of Corollaries \ref{Co1} - \ref{Co4}}
Suppose the prime factorization of a positive integer $n$ is 
\begin{equation*}
	n = \prod_{p|n} p^{\nu_p(n)}.
\end{equation*}
Then, the sum of all positive divisors of $n$ is given by
\begin{equation*}
\sigma(n) = \prod_{p|n} (1+p+p^2+p^3+\cdots+p^{\nu_p(n)}).
\end{equation*}
If $p\equiv -1 \pmod{2^j}$ and $\nu_p(n) \equiv -1 \pmod{2^{4-j}}$, where $1\leq j <4$, then we can easily see that
\begin{equation*}
	1+p+p^2+\cdots+p^{\nu_p(n)} \equiv 0 \pmod{2^3}.
\end{equation*}
However, if $p\equiv -1 \pmod{2^j}$ and $\nu_p(n) \equiv -1 \pmod{2^{5-j}}$, where $1\leq j <5$, then
\begin{equation*}
	1+p+p^2+\cdots+p^{\nu_p(n)} \equiv 0 \pmod{2^4}.
\end{equation*}
Corollaries \eqref{Co1} and \eqref{Co2} follow from Theorem \ref{T2} and the above arguments. The immediate consequences of these corollaries are corollaries \ref{Co3} and \ref{Co4}.

\end{document}